\magnification=1202

\input amssym.def
 at9.98pt
\font\bb=msbm10 at9.98pt
 at9.98pt 
 at5pt
\font\cyr=wncyi10 at9.98pt
\font\eightrm=cmr8
\font\eightsc=cmcsc8
\font\labf=cmbx10 at13.1pt
 at15.74pt
\font\larm=cmr10 at13.1pt
 at15.74pt
\font\sc=cmcsc10 at9.98pt
\font\tenpbf=cmbx10 at8.32pt
\font\tenpit=cmti10 at8.32pt
\font\tenprm=cmr10 at8.32pt


\def\Ad{{\rm Ad}}

\def\an{\raise0.5pt\hbox{$\kern2pt\scriptstyle\in\kern2pt$}}
\def\Ann{\hbox{\rm Ann\kern1pt}}
\def\Anns{\hbox{$\scriptstyle\rm Ann\kern0.5pt$}}

\def\arkef{\advance\chapternumber by 1\sc\roman{\the\chapternumber}}
\def\Aut{\hbox{\rm Aut\kern1pt}}
\def\bell{\hskip0pt\lower1.6pt\hbox{\bel\char'012}\kern5pt}

\def\callige{\hbox{\calligl e\kern2pt}}
\def\cheridexi{\hskip0pt\lower2pt\hbox{\cheridexia}\kern5pt}
\def\cheridexia{{\bbding\char'21}}

\def\Coad{{\rm Coad}}
\def\coker{{\rm coker\kern1pt}}
\def\Colon{\colon\kern2pt}
\def\comp{\hbox{\lower5.8pt\hbox{\larm\char'027}}}

\def\corang{{\rm corang\kern1pt}}
\def\cos{\hbox{\rm cos\kern1pt}}
\def\cosh{\hbox{\rm cosh\kern1pt}}
\def\dbaraux{\hbox{\= {\kern-2pt\= {}}}}
\def\dbar#1{\raise3pt\hbox{\dbaraux}\kern-7.8pt #1}
\def\Der{\lower0.5pt\hbox{\ygoth Der}}

\def\dim{{\rm dim\kern1pt}}
\def\double{\hbox{\kern1.5pt\bb\char'156\kern-7.6pt\char'157\kern1.5pt}}

\def\enwsh#1{{\lower2.1pt\hbox{$\buildrel{\textstyle\cup}\over
{\lower.8pt\hbox{${}_{\scriptscriptstyle#1}$}}$}}}
\def\exp{{\rm exp\kern1pt}}

\def\exten{\hbox{\callig \kern-2.5pt Ext\lower2.5pt\hbox{\kern2.5pt}}}
\def\Ham{\hbox{\rm Ham\kern1pt}}
\def\im{{\rm im\kern1pt}}
\def\k{\raise0.25pt\hbox{$\ygot k$}}

\def\ker{{\rm ker\kern1pt}}

\def\Lie{\hbox{
\callig Lie\kern2pt}}
\def\mavrodexi{\hskip0pt\lower2pt\hbox{\mavrodexia}\kern5pt}
\def\mavrodexia{{\bbding\char'15}}
\def\meriki{\hbox{\cyr\char'144\kern0.3pt}}

\def\na{\raise0.5pt\hbox{$\kern2pt\scriptstyle\ni\kern2pt$}}
\def\noan{\hbox{$\an\raise0.6pt\hbox{$\kern-6.5pt\scriptstyle
          \slash\kern3pt$}$}}

\def\oplus{\;{\mathchar"2208}\;}

\def\otimes{\;{\mathchar"220A}\;}
\def\pounds{\rlap{\lower3.5pt\hbox{\kern2.9pt\hbox{\char'26}}}
           {\script L}}
\def\pr{\hbox{\kern3pt{\calligs p}\callig r\kern2pt}}
\def\qed{\hbox{\kern0.3cm\vrule height5pt width5pt depth-0.2pt}}
\def\QED{\hbox{\kern0.3cm\vrule height6pt width6pt depth-0.2pt}}

\def\rang{{\rm rang\kern1pt}}
\def\rank{{\rm rank\kern1pt}}

\def\san{\raise0.5pt\hbox{$\kern0.7pt\scriptscriptstyle
         \in\kern0.7pt$}}

\def\scomp{\hskip-0.05truecm\hbox{\lower5pt\hbox{$\mathchar"2017$}}
           \hskip-0.05truecm}

\def\sem{\hbox{{\script S}\kern-2.5pt\callig em\kern2pt}}
\def\semidir{\hbox{$\;$\bb\char'156$\;$}}
\def\sin{\hbox{\rm sin\kern1pt}}
\def\sinh{\hbox{\rm sinh\kern1pt}}

\def\styl{\hbox{\bbding\char'26}}
\def\stylo{\hskip0.3truecm\hbox{\lower1.5pt\hbox{\styl}}}
\def\times{\;{\mathchar"2202}\;}

\def\tonos{\hbox{\kern-1.3pt\lower0.7pt\hbox{$\mathchar"6013$}}}
\def\tonoskef{\hbox{$\kern-1.3pt\mathchar"6013$}}
\def\wedget{\hbox{$\;{\scriptstyle\wedge}\;$}}
\def\wbaraux{\hbox{\= {\kern-1.4pt\= {\kern-1.4pt\= {\kern-1.4pt\=
 {\kern-1.4pt\= {\kern-1.4pt\= {\kern-1.4pt\= {\kern-1.4pt\= {}}}}}}}}}}
\def\wbar#1{\hbox{\raise3pt\hbox{\wbaraux}\kern-30.5pt #1}}

\def\wwbaraux{\hbox{\= {\kern-1.4pt\= {\kern-1.4pt\= {\kern-1.4pt\=
{\kern-1.4pt\= {\kern-1.4pt\= {\kern-1.4pt\= {\kern-1.4pt\=
{\kern-1.4pt\= {}}}}}}}}}}}
\def\wwbar#1{\hbox{\raise3pt\hbox{\wwbaraux}\kern-34pt #1}}

\catcode`\@=11
\def\eightpoint{\eightrm}
\def\footnote#1{\edef\@sf{\spacefactor\the\spacefactor}#1\@sf
     \insert\footins\bgroup\eightpoint
     \interlinepenalty100 \let\par=\endgraf
      \leftskip=0pt \rightskip=0pt
      \splittopskip=10pt plus 1pt minus 1pt \floatingpenalty=20000
      \smallskip\item{#1}\bgroup\strut\aftergroup\@foot\let\neft}
\skip\footins=12pt plus 2pt minus 4pt
\dimen\footins=30pc

\def\line{\hbox to\hsize}

\def\title#1{\line{\hss}\line{\hss#1\hss}%
\line{\hss}\hskip-0.75truecm}

\def\author#1{{\tenprm #1:}}
\def\ekdoths#1{{\tenprm #1}}

\def\periodiko#1{{\tenpit #1\tenprm ,}}
\def\selides#1{{\tenprm #1}}
\def\titlosa#1{{\tenprm #1,}}
\def\titlosb#1{{\tenpit #1\tenprm ,}}
\def\volume#1{{\tenprm Vol. \tenpbf #1\tenprm :}}

%
%
%
\def\teleia{\hbox{.}}
\newif\ifPhysRev
\def\Textindent#1{\noindent\llap{#1\enspace}\ignorespaces}
\def\nonfrenchspacing{\sfcode`\.=3001 \sfcode`\!=3000 \sfcode`\?=3000
        \sfcode`\:=2000 \sfcode`\;=1500 \sfcode`\,=1251 }
\nonfrenchspacing
\newdimen\d@twidth
 {\setbox0=\hbox{s.} \global\d@twidth=\wd0 \setbox0=\hbox{s}
        \global\advance\d@twidth by -\wd0 }
\def\removehglue{\loop \unskip \ifdim\lastskip >\z@ \repeat }
\def\roll@ver#1{\removehglue \nobreak \count255 =\spacefactor \dimen@=\z@
        \ifnum\count255 =3001 \dimen@=\d@twidth \fi
        \ifnum\count255 =1251 \dimen@=\d@twidth \fi
    \iftwelv@ \kern-\dimen@ \else \kern-0.83\dimen@ \fi
   #1\spacefactor=\count255 }
\def\step@ver#1{\relax \ifmmode #1\else \ifhmode
        \roll@ver{${}#1$}\else {\setbox0=\hbox{${}#1$}}\fi\fi }
\def\attach#1{\step@ver{\strut^{\mkern 2mu #1} }}

\normalbaselineskip = 20pt plus 0.2pt minus 0.1pt
\normallineskip = 1.5pt plus 0.1pt minus 0.1pt
\normallineskiplimit = 1.5pt
\newskip\normaldisplayskip
\normaldisplayskip = 20pt plus 5pt minus 10pt
\newskip\normaldispshortskip
\normaldispshortskip = 6pt plus 5pt
\newskip\normalparskip
\normalparskip = 6pt plus 2pt minus 1pt
\newskip\skipregister
\skipregister = 5pt plus 2pt minus 1.5pt
\newif\ifsingl@    \newif\ifdoubl@
\newif\iftwelv@    \twelv@true
\def\singlespace{\singl@true\doubl@false\spaces@t}
\def\doublespace{\singl@false\doubl@true\spaces@t}
\def\normalspace{\singl@false\doubl@false\spaces@t}
\def\Tenpoint{\tenpoint\twelv@false\spaces@t}
\def\Twelvepoint{\twelvepoint\twelv@true\spaces@t}
\def\spaces@t{\relax
      \iftwelv@ \ifsingl@\subspaces@t3:4;\else\subspaces@t1:1;\fi
       \else \ifsingl@\subspaces@t3:5;\else\subspaces@t4:5;\fi \fi
      \ifdoubl@ \multiply\baselineskip by 5
         \divide\baselineskip by 4 \fi }
\def\subspaces@t#1:#2;{
      \baselineskip = \normalbaselineskip
      \multiply\baselineskip by #1 \divide\baselineskip by #2
      \lineskip = \normallineskip
      \multiply\lineskip by #1 \divide\lineskip by #2
      \lineskiplimit = \normallineskiplimit
      \multiply\lineskiplimit by #1 \divide\lineskiplimit by #2
      \parskip = \normalparskip
      \multiply\parskip by #1 \divide\parskip by #2
      \abovedisplayskip = \normaldisplayskip
      \multiply\abovedisplayskip by #1 \divide\abovedisplayskip by #2
      \belowdisplayskip = \abovedisplayskip
      \abovedisplayshortskip = \normaldispshortskip
      \multiply\abovedisplayshortskip by #1
        \divide\abovedisplayshortskip by #2
      \belowdisplayshortskip = \abovedisplayshortskip
      \advance\belowdisplayshortskip by \belowdisplayskip
      \divide\belowdisplayshortskip by 2
      \smallskipamount = \skipregister
      \multiply\smallskipamount by #1 \divide\smallskipamount by #2
      \medskipamount = \smallskipamount \multiply\medskipamount by 2
      \bigskipamount = \smallskipamount \multiply\bigskipamount by 4 }
\def\normalbaselines{ \baselineskip=\normalbaselineskip
   \lineskip=\normallineskip \lineskiplimit=\normallineskip
   \iftwelv@\else \multiply\baselineskip by 4 \divide\baselineskip by 5
     \multiply\lineskiplimit by 4 \divide\lineskiplimit by 5
     \multiply\lineskip by 4 \divide\lineskip by 5 \fi }


\def\abstract#1{\parshape=1 0.7cm \dimen10
                {\tenpbf Abstract. \tenprm #1}}

\newcount\appendixnumber     \appendixnumber=0
\newcount\chapternumber      \chapternumber=0
\newcount\equanumber         \equanumber=0
\newcount\mathnumber         \mathnumber=0
\newcount\appequanumber      \appequanumber=0
\newcount\appmathnumber      \appmathnumber=0

\let\variableone=\relax
\let\variabletwo=\relax
\let\chapterlabel=\relax
\let\sectionlabel=\relax
\let\mathlabel=\relax
\newtoks\chapterstyle        \chapterstyle={\Number}
\newtoks\sectionstyle        \sectionstyle={\chapterlabel\Number}
\newskip\chapterskip         \chapterskip=\bigskipamount
\newskip\sectionskip         \sectionskip=\medskipamount
\newskip\headskip            \headskip=8pt plus 3pt minus 3pt
\newdimen\chapterminspace    \chapterminspace=15pc
\newdimen\sectionminspace    \sectionminspace=10pc
\newdimen\sectionspace       \sectionspace=20pc
\newdimen\referenceminspace  \referenceminspace=25pc

\def\chapterreset{\global\advance\chapternumber by 1
   \ifnum\equanumber<0 \else\global\equanumber=0\fi
   \mathnumber=0
   \makechapterlabel}
\def\makechapterlabel{\let\sectionlabel=\relax\let\mathlabel=\relax
 \xdef\chapterlabel{\the\chapterstyle{\the\chapternumber\teleia\kern3pt}}}

\def\rightheadline{\sc\hfil\variableone\eightsc\hfil\folio}
\def\leftheadline{\eightsc\folio\hfil{\sc\variabletwo}\hfil}
\def\heads{\footline={\hfil}\headline={\ifodd\pageno
               \rightheadline\else\leftheadline\fi}}

\def\headseis{\partreset\headline={\ifodd\pageno{
                         \hfil\sc partie {\eightsc\the\partnumber}
                         -introduction\hfil\eightsc\folio}\else
                        {\eightsc\folio\hfil\sc partie
                         {\eightsc\the\partnumber}-introduction\hfil}\fi}
                        \footline={\hfil}}

\def\alphabetic#1{\count255='140 \advance\count255 by #1\char\count255}
\def\Alphabetic#1{\count255='100 \advance\count255 by #1\char\count255}
\def\Roman#1{\uppercase\expandafter{\romannumeral #1}}
\def\roman#1{\romannumeral #1}
\def\Number#1{\number #1}
\def\BLANC#1{}

\def\titlestyle#1{\par\begingroup \interlinepenalty=9999
     \leftskip=0.02\hsize plus 0.23\hsize minus 0.02\hsize
     \rightskip=\leftskip \parfillskip=0pt
     \hyphenpenalty=9000 \exhyphenpenalty=9000
     \tolerance=9999 \pretolerance=9000
     \spaceskip=0.333em \xspaceskip=0.5em
     \iftwelv@\bf\else\bf\fi
   \noindent #1\par\endgroup }

\def\spacecheck#1{\dimen@=\pagegoal\advance\dimen@ by -\pagetotal
   \ifdim\dimen@<#1 \ifdim\dimen@>0pt \vfil\break \fi\fi}
\def\TableOfContentEntry#1#2#3{\relax}

\def\chapter#1{\par\vskip0.7cm
   \chapterreset \titlestyle{\chapterlabel\ #1}
   \nobreak\vskip\headskip
   \wlog{\string\chapter\space \chapterlabel} }

\def\appendixreset{\global\advance\appendixnumber by 1
                   \appmathnumber=0\appequanumber=0}
\def\appendix#1{\par \penalty-300\vskip\chapterskip
   \spacecheck\chapterminspace
   \appendixreset \title{\bf Appendix \Alphabetic{\the\appendixnumber}}
   \nobreak\vskip-\chapterskip\penalty 30000
   \vskip-\chapterskip
   \par{\titlestyle{#1}}
   \vskip\chapterskip
   \wlog{\string\appendix\space \chapterlabel} }

%
%
\def\eqname#1{\relax \ifnum\equanumber<0
     \xdef#1{{\noexpand\rm(\number-\equanumber)}}%
       \global\advance\equanumber by -1
    \else \global\advance\equanumber by 1
      \xdef#1{{\noexpand(\rm{\the\chapternumber}\teleia
                            \rm{\number\equanumber})}} \fi #1}

\def\eqn{\eqno\eqname}

\def\math#1#2{\vskip0.1cm
   \global\advance\mathnumber by 1
   \xdef\mathlabel{\the\chapternumber\teleia\the\mathnumber}
   \wlog{\string\math\space \mathlabel}
   {\bf\enspace\mathlabel\hskip0.2cm #1}
   \xdef#2{{\mathlabel}}}

\def\appeqname#1{\relax \ifnum\appequanumber<0
     \xdef#1{{\noexpand\rm(\number-\appequanumber)}}%
       \global\advance\appequanumber by -1
    \else \global\advance\appequanumber by 1
      \xdef#1{{\noexpand(\hbox{\Alphabetic{\the\appendixnumber}}\teleia
                            {\number\appequanumber})}} \fi #1}

\def\mathapp#1#2{\vskip0.1cm
   \global\advance\appmathnumber by 1
   \xdef\appmathlabel{{\Alphabetic{\the\appendixnumber}}\teleia
   \the\appmathnumber}
   \wlog{\string\mathapp\space \appmathlabel}
   {\bf\enspace\appmathlabel\hskip0.2cm #1}
   \xdef#2{{\appmathlabel}}}


%
%
%
\newtoks\referencestyle      \referencestyle={\tenpbf\Number}
\newcount\referencecount     \referencecount=0
\newcount\lastrefsbegincount \lastrefsbegincount=0
\newif\ifreferenceopen       \newwrite\referencewrite
\newif\ifrw@trailer
\newdimen\refindent     \refindent=13pt
\def\NPrefmark#1{\attach{\scriptscriptstyle [ #1 ] }}
\let\PRrefmark=\attach
\def\refmark#1{\relax\ifPhysRev\PRrefmark{#1}\else\NPrefmark{#1}\fi}
\def\refend@{\refmark{\number\referencecount}}
\def\refend{\refend@{}\space }
\def\refsend{\refmark{\count255=\referencecount
   \advance\count255 by-\lastrefsbegincount
   \ifcase\count255 \number\referencecount
   \or \number\lastrefsbegincount,\number\referencecount
   \else \number\lastrefsbegincount-\number\referencecount \fi}\space }
\def\refitem#1{\par\hangafter=0 \hangindent=\refindent	\Textindent{#1}}
\def\Ref{\rw@trailertrue\REF}
\def\REF#1{\r@fstart{#1}%
   \rw@begin{\tenprm [\tenpbf\Number{\the\referencecount}\tenprm ]}\rw@end}
\def\r@fstart#1{\chardef\rw@write=\referencewrite \let\rw@ending=\refend@
   \ifreferenceopen \else \global\referenceopentrue
   \immediate\openout\referencewrite=referenc.txa
   \toks0={\catcode`\^^M=10}\immediate\write\rw@write{\the\toks0} \fi
   \global\advance\referencecount by 1 
   \xdef#1{[{\the\referencestyle{\the\referencecount}}]}}
 {\catcode`\^^M=\active %
 \gdef\rw@begin#1{\immediate\write\rw@write{\noexpand\refitem{#1}}%
   \begingroup \catcode`\^^M=\active \let^^M=\relax}%
 \gdef\rw@end#1{\rw@@end #1^^M\rw@terminate \endgroup%
   \ifrw@trailer\rw@ending\global\rw@trailerfalse\fi }%
 \gdef\rw@@end#1^^M{\toks0={#1}\immediate\write\rw@write{\the\toks0}%
   \futurelet\n@xt\rw@test}%
 \gdef\rw@test{\ifx\n@xt\rw@terminate \let\n@xt=\relax%
       \else \let\n@xt=\rw@@end \fi \n@xt}%
}
\let\rw@ending=\relax
\let\rw@terminate=\relax

\def\Closeout#1{\toks0={\catcode`\^^M=5}\immediate\write#1{\the\toks0}%
   \immediate\closeout#1}
\topskip1truecm
\voffset=2.5truecm
\hsize 15truecm
\vsize 20truecm
\hoffset=0.5truecm
\def\undertext#1{$\underline{\hbox{#1}}$}
\topglue 3truecm
\nopagenumbers
\def\variableone{p. baguis}
\def\variabletwo{induction of hamiltonian poisson actions}
\def\absize{11cm}

\def\abstract#1{\baselineskip=12pt plus .2pt
                \parshape=1 0.7cm \absize
                 {\tenpbf Abstract. \tenprm #1}}



\Ref\abmar{
\author{Abraham, R., Marsden, J. E.} 
\titlosb{Foundations of Mechanics}
\ekdoths{Addison-Wesley Publishing Company, Inc. (1978)}}

\Ref\babelon{
\author{Babelon, O., Bernard, D.}
\titlosa{Dressing Symmetries}
\periodiko{Comm. Math. Phys.}
\volume{149}
\selides{279--306 (1992)}}

\Ref\bag{
\author{Baguis, P.}
\titlosa{Semidirect products and the Pukanszky condition}
\periodiko{J. Geom. Phys.}
\volume{25}
\selides{245--270 (1998)}}

\Ref\phd{
\author{Baguis, P.}
\titlosb{Proc\'edures de r\'eduction et d'induction en g\'eom\'etrie
symplectique et de Poisson. Applications.}
\ekdoths{Th\`ese de Doctorat, Universit\'e d'Aix-Marseille II,
D\'ecembre 1997}}

\Ref\drone{
\author{Drinfel'd, V. G.}
\titlosa{Hamiltonian structures on Lie groups, Lie bialgebras
and the geometric meaning of the classical Yang-Baxter equations}
\periodiko{Soviet Math. Dokl.}
\volume{27(1)}
\selides{68--71 (1983)}}

\Ref\drtwo{
\author{Drinfel'd, V. G.}
\titlosa{Quantum groups}
\periodiko{Proc. ICM. Berkeley}
\volume{1}
\selides{789--820 (1986)}}

\Ref\de{
\author{Duval, C., Elhadad, J.}
\titlosa{Geometric quantization and localization of
relativistic spin systems}
\periodiko{Contemp. Math.}
\volume{132}
\selides{317--330 (1992)}}

\Ref\det{
\author{Duval, C., Elhadad, J., Tuynman, G. M.}
\titlosa{Pukanszky's condition and symplectic induction}
\periodiko{J. Diff. Geom.}
\volume{36}
\selides{331--348 (1992)}}

\Ref\gelone{
\author{Gel'fand, I. M., Dorfman, I. Ya.}
\titlosa{Hamiltonian operators and Algebraic structures related to 
them}
\periodiko{Funct. Anal. Appl.}
\volume{13}
\selides{248--262 (1979)}}

\Ref\geltwo{
\author{Gel'fand, I. M., Dorfman, I. Ya.}
\titlosa{The Schouten bracket and Hamiltonian operators}
\periodiko{Funct. Anal. Appl.}
\volume{14}
\selides{223--226 (1980)}}

\Ref\kazhdan{
\author{Kazhdan, D., Kostant, B., Sternberg, S.}
\titlosa{Hamiltonian Group Actions and Dynamical Systems
of Calogero Type}
\periodiko{Comm. Pure Appl. Math.}
\volume{31}
\selides{481--508 (1978)}}

\Ref\konts{
\author{Kontsevich, M.}
\titlosa{Deformation quantization of Poisson manifolds, I}
\periodiko{q-alg/9709040}
\selides{(1997)}}

\Ref\jhlu{
\author{Lu, J.-H.} 
\titlosb{Multiplicative and affine
Poisson structures on Lie Groups} 
\ekdoths{Ph.D. thesis, Univ. of California, Berkeley (1990)}}

\Ref\luwei{
\author{Lu, J.-H., Weinstein, A.}
\titlosa{Poisson-Lie groups, dressing transformations and Bruhat
decompositions}
\periodiko{J. Diff. Geom.}
\volume{31}
\selides{501--526 (1990)}}

\Ref\magri{
\author{Magri, F.}
\titlosa{A simple model of the integrable Hamiltonian equation}
\periodiko{J. Math. Phys.}
\volume{19(5)}
\selides{1156--1162 (1978)}}

\Ref\mamo{
\author{Magri, F., Morosi, C.}
\titlosa{A geometric characterization of integrable Hamiltonian 
systems
through the theory of Poisson-Nijenhuis manifolds}
\periodiko{Quaderno S 19}
\selides{1978, Universit\`a di Milano}}

\Ref\mr{
\author{Marsden, J. E., Ratiu, T.}
\titlosa{Reduction of Poisson manifolds}
\periodiko{Lett. Math. Phys.}
\volume{11}
\selides{161--169 (1986)}}

\Ref\marwei{
\author{Marsden, J. E., Weinstein, A.}
\titlosa{Reduction of symplectic manifolds with symmetry}
\periodiko{Rep. Math. Phys.}
\volume{5}
\selides{121--130 (1974)}}

\Ref\petalone{
\author{Petalidou, F.}
\titlosb{\'Etude locale de structures bihamiltoniennes}
\ekdoths{Th\`ese de Doctorat, Universit\'e Pierre et Marie 
Curie-Paris VI,
1998}}

\Ref\semenov{
\author{Semenov-Tian-Shansky, M. A.}
\titlosa{Dressing transformations and Poisson-Lie group
actions}
\periodiko{Publ. RIMS, Kyoto University}
\volume{21}
\selides{1237--1260 (1985)}}

\Ref\souriau{
\author{Souriau, J.-M.}
\titlosb{Structures des syst\`emes dynamiques}
\ekdoths{Dunod, Paris (1969)}}

\Ref\va{
\author{Vaisman, I.} 
\titlosb{Lectures on the Geometry of
Poisson Manifolds}
\ekdoths{Birkh\"auser Verlag (1994) {\tenpit and cited references}}}

\Ref\xu{
\author{Xu, P.}
\titlosa{Symplectic groupoids of reduced Poisson spaces}
\periodiko{C. R. Acad. Sci. Paris}
\volume{314}
\selides{457--461, S\'erie I (1992)}}


\hskip10cm

\centerline{\labf Induction of Hamiltonian Poisson actions}

\vskip1cm

\centerline{\bf P. Baguis\footnote{$^{1}$}{e-mail: 
pbaguis@ulb.ac.be}}

\vskip0.3cm

\centerline{Universit\'e Libre de Bruxelles}
\centerline{Campus Plaine, CP 218 Bd du Triomphe}
\centerline{1050, Brussels, Belgium}

\vskip1cm

\abstract{
We propose a Poisson-Lie analog of the symplectic 
induction procedure, using an appropriate Poisson generalization of
the reduction of symplectic manifolds with symmetry. 
Having as basic tools the equivariant momentum maps of Poisson 
actions, the double group of a Poisson-Lie group and the reduction of Poisson
manifolds with symmetry,  we show how
one can induce a Poisson action admitting an equivariant momentum
map. We prove that, under certain conditions,  the dressing orbits of  
a Poisson-Lie group can be obtained by Poisson induction from the 
dressing orbits of a Poisson-Lie subgroup.}

\vskip1cm

{\tenprm 
{\tenpit Key-words}: Poisson-Lie groups, induction of Poisson 
actions, dressing orbits

\vskip0.2cm

{\tenpit 1991 MSC}: 53C15}

\vfill\eject

\baselineskip=14pt plus .2pt

\chapter{Introduction}

Poisson manifolds occur as phase spaces in Hamiltonian mechanics and 
have
important applications to the theory of completely integrable systems.
This is, in particular, the case of bihamiltonian manifolds, that is 
manifolds
equipped with two Poisson structures $\pi_{1}$ and $\pi_{2}$ such that 
$[\pi_{1},\pi_{2}]=0$, see {\gelone}, {\geltwo}, {\magri}, {\mamo},
{\petalone}. The algebras of
observables in quantum mechanics are also relevant to Poisson 
geometry, as
explained in {\konts}.

A Lie group equipped with a Poisson structure such that the 
corresponding 
group operation be a Poisson map, is called Poisson-Lie group. This
particularly interesting and rich structure has first been studied in
{\drone} and {\semenov} (see also {\luwei} and the monograph {\va}).
Poisson-Lie groups arise naturally in problems of quantum field theory
and integrable systems. For example, a solution of the quantum
Yang-Baxter equation defines a ``quantum group'' in the sense of 
{\drtwo}
which, by definition, is a Hopf algebra. Formally, the ``classical 
limit''
of a quantum group is a Poisson-Lie group. 

\heads

On the other hand, there exist integrable systems, as for example the
KdV equations, for which Poisson-Lie groups provide a deeper insight.
For such systems, the dressing tansformation groups play the r\^ole
of ``hidden symmetry'' groups. According to {\semenov}, the dressing
transformation group does not in general preserve the Poisson 
structure
on the phase space. Furthermore, it carries a natural Poisson 
structure
defined by the Riemann-Hilbert problem entering the definition of the
dressing transformations, and it turns out that this Poisson structure
makes the dressing transformation group into a Poisson-Lie group.

In the same context, the Hamiltonian actions of Poisson-Lie groups
have clarified several aspects of the soliton equations. Indeed, the
dressing transformations of the soliton equations which admit a Lax
representation, are generated by the monodromy matrix {\babelon}, 
which 
in this case is a momentum mapping in the sense of {\jhlu}.

Our aim in this paper is to generalize and study the procedure of
symplectic induction {\kazhdan}, {\de}, {\det} 
in the context of Poisson-Lie groups and
Poisson manifolds. As we shall explain, this generalization is
possible in the following sense: given a Poisson-Lie group 
$(G,\pi_{G})$,
a Poisson-Lie subgroup $(H,\pi_{H})\hookrightarrow(G,\pi_{G})$, a 
Poisson
manifold $(P,\pi_{P})$ and a Hamiltonian action $H\times P\rightarrow 
P$
with equivariant momentum mapping $P\rightarrow H^{\ast}$, one can 
construct a new Poisson manifold $(P_{ind},\pi_{ind})$ on which the
Poisson-Lie group $(G,\pi_{G})$ acts in a Hamiltonian way. This
statement is our basic result and it is given by Theorem
4.3.  As in the symplectic case, an appropriate reduction procedure 
(for Poisson manifolds now) is needed. This is easily obtained putting 
together known facts about Poisson reduction {\jhlu}, {\va}, see 
Theorem 2.1. We also need appropriate Poisson generalizations of the 
natural Hamiltonian actions of a Lie group $G$ and a Lie subgroup 
$H\subset G$ on the cotangent bundle $T^{\ast}G$ from which the 
induced manifold is constructed {\det}. Propositions 3.1 and 3.3 
describe these actions in the Poisson setting.

We finally prove that the Poisson induction procedure can be used in 
order to find Poisson generalizations of the modified cotangent 
bundles {\det} and of the symplectic induction of coadjoint orbits 
{\bag} {\phd}.

\vskip0.2cm

{\bf Conventions.} If $(P,\pi_{P})$ is a Poisson manifold, then 
$\pi_{P}^{\sharp}\colon T^{\ast}P\rightarrow TP$ is the map defined 
by 
$\alpha(\pi_{P}^{\sharp}(\beta))=\pi_{P}(\alpha,\beta),\forall\alpha,\beta\an 
T^{\ast}P$. Let now $\sigma\colon G\times 
P\rightarrow P$ (resp.  $\sigma\colon P\times 
G\rightarrow P$) be a left (resp. right) Poisson action of the Poisson-Lie group
$(G,\pi_{G})$ on $(P,\pi_{P})$,  and let us denote by $\sigma(X)$ the 
infinitesimal generator of the action and by $G^{\ast}$ the dual group 
of $G$. Then, we say that 
$\sigma$ is Hamiltonian if there exists a differentiable map $J\colon 
P\rightarrow G^{\ast}$,  called momentum mapping, 
satisfying the following equation, for each $X\an\frak g$:
$$\sigma(X)=\pi_{P}^{\sharp}(J^{\ast}X^{l})\quad(resp. \quad
\sigma(X)=-\pi_{P}^{\sharp}(J^{\ast}X^{r})).$$
In the previous equation $X^{l}$ (resp. $X^{r}$) is the left (resp. 
right) invariant 1-form on $G^{\ast}$ whose value at the identity is 
equal to $X\an\frak g\cong(\frak g^{\ast})^{\ast}$. The momentum 
mapping is said to be equivariant, if it is a morphism of Poisson 
manifolds with respect to the Poisson structure $\pi_{P}$ on $P$ and 
the canonical Poisson structure on the dual group of the Poisson Lie 
group $(G,\pi_{G})$. Left and right 
infinitesimal dressing actions $\lambda\colon\frak 
g^{\ast}\rightarrow{\cal X}(G)$ and $\rho\colon\frak 
g^{\ast}\rightarrow{\cal X}(G)$ of $\frak g$ on $G^{\ast}$ are defined by
$$\lambda(\xi)=\pi_{G}^{\sharp}(\xi^{l})\quad\hbox{and}\quad
\rho(\xi)=-\pi_{G}^{\sharp}(\xi^{r}),\quad\forall \xi\an\frak g^{\ast}.$$
Similarly, one defines infinitesimal left and right dressing actions 
of $\frak g$ on $G^{\ast}$. In the case where the vector fields 
$\lambda(\xi)$ (or, equivalently, $\rho(\xi)$) are complete for all 
$\xi\an\frak g^{\ast}$,  we have left and right actions of 
$(G^{\ast},\pi_{G^{\ast}})$ on $(G,\pi_{G})$ denoted also by $\lambda$ 
and $\rho$ respectively, and we say that $(G,\pi_{G})$
is a complete Poisson-Lie group.

\vskip0.2cm

{\bf Acknowledgments.} It is a pleasure to thank Professor C. Duval 
for
careful and critical reading of the manuscript
and for many stimulating discussions during the preparation of this 
work. I would like also to thank Professor M. Cahen for his interest 
in this work and for his comments on the manuscript.

\chapter{Reduction of Poisson manifolds}

The reduction of symplectic manifolds with symmetry has been
systematically studied in {\marwei}. The importance 
of this procedure for Hamiltonian dynamics is already
very clear as it describes in a unified way several properties
of Hamiltonian systems. The Poisson generalization of reduction
with symmetry has been carried out in {\jhlu} for the special case
of a Poisson action of a Poisson-Lie group on a symplectic
manifold, admiting a momentum map. On the other hand, reduction 
of Poisson manifolds with symmetry under the Hamiltonian action of an 
ordinary Lie group can be found in {\va}. Here we will study a somewhat more
general situation where a Poisson-Lie group acts in a 
Hamiltonian way on a Poisson manifold. Before we state
the reduction theorem for Poisson manifolds with symmetry,
we recall the notion of sub-characteristic distribution. If 
$(P,\pi_{P})$ is a Poisson manifold and $N$ a submanifold of $P$,
then we define the sub-characteristic distribution of $N$ as
$${\cal C}N=\pi_{P}^{\sharp}((TN)^{\circ})\cap TN,\eqn\subchar$$
where $(TN)^{\circ}$ is the annihilator of the tangent bundle $TN$: 
$$(T_{x}N)^{\circ}=\{\alpha\an T_{x}^{\ast}P
\;|\;\alpha(v)=0,\;\forall v\an T_{x}N\}.$$
We will deal only with Poisson actions of Poisson-Lie groups admiting 
equivariant momentum mappings. Although this seems to be a strong 
condition on the Poisson action, it has been proved {\phd} that, at 
least for Poisson actions on symplectic manifolds, one can, under 
reasonable conditions, be reduced to the equivariant case.

\math{Theorem.}{\reductionth}{\sl Let $(P,\pi_{P})$ be a Poisson
manifold and $\sigma\colon G\times P\rightarrow P$ a 
Poisson action of the connected Poisson-Lie group $(G,\pi_{G})$ on 
$(P,\pi_{P})$ admiting an equivariant momentum mapping 
$J\colon P\rightarrow G^{\ast}$. Let $u\an G^{\ast}$ be an element
such that: (1) $u$ is a regular value for all the restrictions of $J$
to the symplectic leaves of $P$; (2) the submanifold $J^{-1}(u)$
has a clean intersection with the symplectic leaves of $P$. Then,
if $G_{u}$ is the isotropy subgroup of $u$ with respect to the
left dressing action of $G$ on $G^{\ast}$, the sub-characteristic 
distribution of $J^{-1}(u)$ defines a regular foliation (that is 
of constant dimension) whose leaves are the orbits of $G_{u}$. 
Furthermore, if this foliation is defined by a submersion  
$s\colon J^{-1}(u)\rightarrow P_{u}$, then the manifold $P_{u}$
possesses a well-defined Poisson structure whose symplectic 
distribution is the projection of ${\cal S}(P)\cap TJ^{-1}(u)$,
where ${\cal S}(P)$ is the symplectic distribution 
of $(P,\pi_{P})$.}

\undertext{\it Proof}. We observe that the existence of a momentum 
mapping for the action $\sigma$, implies that the orbit $G\cdot x$,
for each $x\an P$, is contained in the symplectic leaf $S(x)$ through $x$ and
for each $x\an J^{-1}(u)$, the orbit $G_{u}\cdot x$ is contained
in $S(x)_{u}=S(x)\cap J^{-1}(u)$. Furthermore, we have 
$\pi_{P}^{\sharp}(x)
\big((T_{x}J^{-1}(u))^{\circ}\big)=T_{x}(G\cdot x)$
and the submanifold $J^{-1}(u)$ 
has a clean intersection with the orbits of $G$ in $P$:
$T(G\cdot x)\cap TJ^{-1}(u)=T(G_{u}\cdot x)$. After these remarks, 
the details of the proof are as in {\jhlu} and {\va}.\QED

\vskip0.3cm

The reduction described in Theorem {\reductionth} 
is called leafwise reduction because the reduced Poisson 
structure is obtained by reducing each symplectic leaf of 
$P$ by the standard procedure of symplectic geometry.

\chapter{Hamiltonian actions on the double Lie group}

Let $G$ be a Lie group and $i\colon H\hookrightarrow G$ a closed
Lie subgroup. We have a right action of $H$ on $G$ given by right
multiplication, $(g,h)\mapsto gh$, $\forall g\an G,h\an H$, and a 
left action of $G$ on itself given by left multiplication,
$(g,g^{\prime})\mapsto gg^{\prime}$, $\forall g,g^{\prime}\an G$.
The cotangent lifts of these two actions are the basis 
of the symplectic induction {\det} and in the left trivialization
$T^{\ast}G\cong G\times\frak g^{\ast}$ they are given 
by the relations:
$$((g,\mu),h)\mapsto (gh,\Coad(h^{-1})\mu),\forall(g,\mu)\an
T^{\ast}G,h\an H\eqn\rsaction$$
$$(g,(g^{\prime},\mu))\mapsto(gg^{\prime},\mu),\forall g\an G,
(g^{\prime},\mu)\an T^{\ast}G.\eqn\lsaction$$
These actions are Hamiltonian and their equivariant momentum
mappings are respectively given by:
$$T^{\ast}G\na(g,\mu)\mapsto -i^{\ast}\mu\an\frak h^{\ast},
\forall(g,\mu)\an T^{\ast}G,\eqn\rsmommap$$
$$T^{\ast}G\na(g,\mu)\mapsto\Coad(g)\mu\an\frak 
g^{\ast},\forall(g,\mu)\an
T^{\ast}G,\eqn\lsmommap$$
where $\frak h$ is the Lie algebra of $H$ and $i^{\ast}\colon
\frak g^{\ast}\rightarrow\frak h^{\ast}$ is the canonical projection.
We will generalize in this section the previous Hamiltonian actions
in the context of Poisson-Lie groups. This generalization will
provide the basis for Poisson induction, as we will see in the sequel.

Let $(G,\pi_{G})$ be a connected, simply connected and complete
Poisson-Lie group and $i\colon(H,\pi_{H})\hookrightarrow(G,\pi_{G})$
a closed Poisson-Lie subgroup. Then, if $D(G)$ is the double group
of $G$, we find, by Proposition II.36 of {\jhlu}, that the right
action 
$r\colon D(G)\times H\rightarrow D(G)$
given by right multiplication,
$$r(d,h)=dh,\forall d\an D(G),h\an H\eqn\rightHaction$$
is a Poisson action for the symplectic structure $\pi_{+}$ on $D(G)$ 
and the Poisson structure $\pi_{H}$ on $H$. 
We recall here that in the case we are
studying the double group $D(G)$ is gobally isomorphic to the product
$G\times G^{\ast}$ with the group law given by the relation
$$(g,u)\cdot(h,v)=(g\rho_{u^{-1}}(h),\lambda_{h^{-1}}(u)v),\;
\forall(g,u),(h,v)\an D.\eqn\gggroup$$ 
Furthermore, there exist two Poisson structures, $\pi_{+}$ (symplectic)
and $\pi_{-}$ (Poisson-Lie) on $D(G)$ given by 
$$\pi_{\pm}(d)={1\over 2}(R_{d}\pi_{0}\pm L_{d}\pi_{0}),$$
where $\pi_{0}\an\Lambda^{2}\frak d$ is the bivector defined by
$\pi_{0}(\xi_{1}+X_{1},\xi_{2}+X_{2})=\xi_{1}(X_{2})-\xi_{2}(X_{1})$,
$\forall\xi_{i}+X_{i}\an\frak d^{\ast},i=1,2$, see {\jhlu} for more 
details. In the defining equation of $\pi_{\pm}$, $L_{d}$ and $R_{d}$ are the 
extensions to multivector 
fields, of left and right multiplication in $D$.

In fact, the right Poisson action given by $\rightHaction$ is
Hamiltonian:

\math{Proposition.}{\ractionham}{\sl The right Poisson action
given by $\rightHaction$ is Hamiltonian with equivariant
momentum mapping $J_{r}\colon D(G)\rightarrow H^{\ast}$ 
which can be taken equal to $$J_{r}=s\comp i^{\ast}\comp p_{2},$$
where $s\colon H^{\ast}\rightarrow H^{\ast}$ is the inversion
on the dual group $H^{\ast}$, 
$i^{\ast}\colon G^{\ast}\rightarrow H^{\ast}$
is the projection of dual groups induced by the inclusion
$i\colon H\hookrightarrow G$, and $p_{2}\colon D(G)\rightarrow
G^{\ast}$ is the projection onto the second factor.}

\undertext{\it Proof}. The infinitesimal generator of the right 
action $\rightHaction$ is given by the relation
$r(Y)(d)=T_{e}L_{d}(i_{\ast}Y,0)$, $\forall
Y\an\frak h$, $d\an D(G)$. 
Setting now $(J_{r}^{\ast}Y^{r})(d)=
(\eta_{1}+Y_{1})\comp T_{g}L_{g^{-1}}
\comp T_{d}R_{u^{-1}}$, where
$d=gu$, $\eta_{1}+Y_{1}\an{\frak d}^{\ast}
\cong\frak g^{\ast}\oplus\frak g$, and
$\pi_{+}^{\sharp}(J_{r}^{\ast}Y^{r})(d)=
(J^{\ast}_{r}Y^{r})_{d}^{\sharp}$, one finds:
$$(J_{r}^{\ast}Y^{r})^{\sharp}_{d}=
(T_{g}R_{u}\comp T_{e}L_{g})
\big[(Y_{1}-T_{g}L_{g^{-1}}
\lambda(\eta_{1})(g))\oplus(T_{u}R_{u^{-1}}
\rho(Y_{1})(u)-\eta_{1})\big].\eqn\jryr$$
The elements $\eta_{1}$ and $Y_{1}$ of the previous expression are 
calculated using the definition of the momentum map $J_{r}$.
One finds $\eta_{1}=0$ and $Y_{1}=-i_{\ast}\Coad(w^{-1})Y$, 
$w=J_{r}(d)=(i^{\ast}u)^{-1}$. We proceed by recalling the following
useful properties of Poisson-Lie groups {\phd}:

\math{Lemma.}{\coadidentity}{\sl 

(1) If $w=(i^{\ast}u)^{-1}$, then
$i_{\ast}\Coad(w^{-1})Y=\Coad(u)i_{\ast}Y$, $\forall u\an G^{\ast}$,
$Y\an\frak h$. 

(2) The left and right dressing vector fields on the Poisson-Lie group
$(G,\pi_{G})$ are related as follows:
$$\rho(\Coad(g)\xi)(g)=-\lambda(\xi)(g),$$
for each $g\an G$, $\xi\an\frak g^{\ast}$.

(3)  If the map $\phi\colon G\times G^{\ast}\rightarrow
D(G)$ given by $\phi(g,u)=gu$ is a global 
diffeomorphism and $X\an\frak g$, $u\an G^{\ast}$, then
$$\Ad_{D(G)}(u)(X\oplus 0)=T_{e}\rho_{u^{-1}}(X)\oplus
(-T_{u}R_{u^{-1}}\lambda(X)(u)),$$
where $\rho_{u}$ is the right dressing transformation 
of $G^{\ast}$ on $G$ and $\lambda(X)$ 
the infinitesimal generator of the left dressing transformation of 
$G$ on $G^{\ast}$.}






Replacing now in $\jryr$ the values of $\eta_{1}$ and 
$Y_{1}$, using Lemma {\coadidentity} and the fact that the
tangent at the identity of the dressing transformations equals to the
coadjoint representation, we find: 
$$\eqalign{-(J_{r}^{\ast}Y^{r})^{\sharp}_{d}&=(T_{g}R_{u}\comp 
T_{e}L_{g})(\Ad_{D(G)}(u)(i_{\ast}Y\oplus 0))\cr
\hfill&=T_{e}L_{d}(i_{\ast}Y\oplus 0),\cr}$$
which proves that $J_{r}$ is indeed an equivariant (because it is a
Poisson morphism) momentum map for the right action $r$.\QED

\vskip0.3cm

Using analogous techniques, one can prove the following:

\math{Proposition.}{\lactionham}{\sl The left action $l\colon
(G,\pi_{G})\times(D(G),\pi_{+})\rightarrow(D(G),\pi_{+})$ 
given by 
$$l_{k}(d)=\lambda_{u}(k\lambda_{u^{-1}}(g))\cdot u=
\lambda_{\rho_{g^{-1}}(u)}(k)g\cdot u,
\forall k\an G,d=gu\an D(G),\eqn\laction$$
is Hamiltonian with equivariant momentum map 
$J_{l}\colon D(G)\rightarrow G^{\ast}$ 
such that 
$$J_{l}(d)=\rho_{g^{-1}}(u),\forall d=gu\an D(G).\eqn\jleft$$}

\vskip0.3cm

In the trivial case where the Poisson structure $\pi_{G}$ is zero,
one has $G^{\ast}=\frak g^{\ast}$ and the dressing 
transformations of $G^{\ast}$ on $G$ are trivial. 
Furthermore, the dressing transformations of $G$ on $G^{\ast}$ reduce
to the coadjoint action of $G$ on $\frak g^{\ast}$ and the group law
on $D(G)=G\times\frak g^{\ast}$ is simply the semi-direct product
structure on $T^{\ast}G=G\semidir\frak g^{\ast}$.
Then, the Hamiltonian actions and their equivariant 
momentum mappings described in Propositions {\ractionham}
and {\lactionham} reduce to the actions and momentum 
mappings given by the relations $\rsaction$,  $\rsmommap$ and $\lsaction$,
$\lsmommap$ respectively, because $\lambda_{u}=id$, $\forall u\an
G^{\ast}$.

\chapter{Induction of Hamiltonian Poisson actions}

We recall first {\phd} for the reader's convenience some properties of the
Hamiltonian actions of Poisson-Lie groups, very useful in what follows.

\goodbreak

\math{Proposition.}{\lrpoissonmoment}{\sl

(1) Let $\sigma\colon P\times G\rightarrow P$ be a right Poisson 
action of the connected, simply connected and complete Poisson-Lie 
group $(G,\pi_{G})$ on the Poisson manifold $(P,\pi_{P})$, 
admiting an equivariant momentum mapping 
$J\colon P\rightarrow G^{\ast}$. Then, the map 
$\tilde{\sigma}\colon G\times P\rightarrow P$ defined as
$$\tilde{\sigma}(g,p)=\sigma(p,[\lambda_{J(p)}(g)]^{-1}),
\;\forall g\an G,p\an P,\eqn\sigmatilde$$ is a left Poisson action.
Furthermore, $J$ is an equivariant momentum map for $\tilde{\sigma}$.

(2) Let $\sigma_{i}\colon G\times P_{i}\rightarrow P_{i}$, $i=1,2$,
be left Poisson actions admiting equivariant momentum mappings
$J_{i}\colon P_{i}\rightarrow G^{\ast}$, where $G$ is as previously.
Then the map 
$\sigma\colon G\times P\rightarrow P$, 
$P=P_{1}\times P_{2}$, defined by
$$\sigma(g,p)=(\sigma_{1}(\lambda_{J_{2}(p_{2})}(g),p_{1}),
\sigma_{2}(g,p_{2})),\;p=(p_{1},p_{2})\an P,\eqn\twosigma$$
is a left Poisson action with respect to the Poisson 
structure $\pi_{P}=\pi_{1}\oplus\pi_{2}$ on $P$.
Furthermore, $J=\tilde{m}\comp
(J_{1}\times J_{2})\colon P\rightarrow G^{\ast}$ is an equivariant 
momentum mapping for $\sigma$, 
where $\tilde{m}\colon G^{\ast}\times
G^{\ast}\rightarrow G^{\ast}$ is the group 
multiplication in $G^{\ast}$.}

\vskip0.3cm

Consider a Poisson-Lie group $(G,\pi_{G})$ and 
let $i\colon(H,\pi_{H})\hookrightarrow(G,\pi_{G})$ be a closed
Poisson-Lie subgroup. In order to simplify the discussion, we 
assume that $(G,\pi_{G})$ is complete, connected and 
simply connected, so the double group $D(G)$ of $G$ 
will be isomorphic to $G\times G^{\ast}$ with the group law 
given by $\gggroup$.

Let $\sigma\colon(H,\pi_{H})\times(P,\pi_{P})\rightarrow(P,\pi_{P})$
be a left Hamiltonian action of $(H,\pi_{H})$ on the Poisson 
manifold $(P,\pi_{P})$, with equivariant momentum mapping 
$J\colon P\rightarrow H^{\ast}$. By Proposition {\lrpoissonmoment},
we have a left Poisson action $\tilde{r}\colon(H,\pi_{H})
\times(D(G),\pi_{+})\rightarrow(D(G),\pi_{+})$ 
canonically associated to the right Poisson action of Proposition 
{\ractionham}, and if 
$$(\check{P},\pi_{\check{P}})=(P,\pi_{P})\times(D(G),\pi_{+}),
\eqn\Pcheck$$
we also have a left Poisson action $\check{\sigma}\colon
(H,\pi_{H})\times(\check{P},\pi_{\check{P}})\rightarrow
(\check{P},\pi_{\check{P}})$ admiting an equivariant momentum
mapping $\check{J}\colon\check{P}\rightarrow H^{\ast}$ given by
$$\check{J}(p,d)=J(p)J_{r}(d),\;\forall(p,d)\an\check{P}.
\eqn\Jcheck$$
Explicitly, the action $\check{\sigma}$ 
is 
given by $$\check{\sigma}_{h}(p,d)=
\left(\sigma(\lambda_{J_{r}(d)}(h),p),
\tilde{r}_{h}(d)\right),\;\forall(p,d)\an\check{P},h\an 
H,\eqn\sigmacheck$$
where $\tilde{r}$ is the left action with
$$\tilde{r}_{h}(d)=d\lambda_{J_{r}(d)}(h)^{-1},\;\forall d\an D(G),
h\an H.\eqn\rtilde$$
We now observe that the 
momentum mapping 
$\check{J}\colon\check{P}\rightarrow H^{\ast}$
is a submersion, so each element of the dual group $H^{\ast}$
is a regular value for $\check{J}$. In particular, if $e^{\ast}$
is the unit of $H^{\ast}$, then $\check{J}^{-1}(e^{\ast})$ is a
submanifold of $\check{P}$. Using the fact that 
$\pi_{+}$ is symplectic, we find that the symplectic leaves
of $(\check{P},\pi_{\check{P}})$ are of the form $S\times D(G)$,
where $S$ is a symplectic leaf of $P$. This means that $e^{\ast}$
is a regular value for all the restrictions of $\check{J}$ 
to the symplectic leaves of $\check{P}$. 

Next, we consider the intersections of the submanifold 
$\check{J}^{-1}(e^{\ast})$ with the symplectic leaves of $\check{P}$.
Using the expression $\Jcheck$ of $\check{J}$ and Proposition
{\ractionham} we find:
$$\check{J}^{-1}(e^{\ast})=\{(p,gu)\an\check{P}=P\times D(G)\;|\;
J(p)=i^{\ast}(u)\}.\eqn\Jcheckminus$$
On the other hand, the symplectic leaf 
$S(m)$ through $m=(p,d)\an\check{P}$ is equal to 
$S(m)=S(p)\times D(G)$, and   
$$\check{J}^{-1}(e^{\ast})\cap S(m)=\{(p,gu)\an S(p)\times D(G)\;|\;
J(p)=i^{\ast}(u)\}.\eqn\Jcheckinter$$
We see now that $T_{n}\check{J}^{-1}(e^{\ast})
\cap T_{n}S(m)=T_{n}(\check{J}^{-1}(e^{\ast})\cap S(m))$, for each
point $n\an\check{J}^{-1}(e^{\ast})\cap S(m)$, which confirms
that $\check{J}^{-1}(e^{\ast})$ has a clean inersection
with the symplectic leaves of $\check{P}$. Furthermore, 
the isotropy subgroup of $e^{\ast}$ with respect to the left
dressing transformations of $H$ on $H^{\ast}$ is the 
group $H$ itself, and if we assume that the action of $H$ on $P$
is proper, then all the conditions of Theorem {\reductionth}
are fulfilled. The quotient manifold 
$$P_{ind}=\check{J}^{-1}(e^{\ast})/H\eqn\Pind$$
which by construction is a Poisson manifold, is called induced 
Poisson manifold. We will denote its Poisson structure as $\pi_{ind}$.

In order to construct a Poisson action of $(G,\pi_{G})$ on
$(P_{ind},\pi_{ind})$, we first study some properties of the
Poisson actions and their momentum mappings of Propositions
{\ractionham} and {\lactionham}. 

\math{Proposition.}{\invar}{\sl 
Let $\check{l}\colon G\times\check{P}
\rightarrow\check{P}$ be the action defined by
$$\check{l}_{k}(p,d)=(p,l_{k}(d)),\;
\forall k\an G,(p,d)\an\check{P},\eqn\lcheckaction$$
where the action $l\colon G\times D(G)\rightarrow D(G)$ 
is given by $\laction$. Then, $\check{l}$ is a Poisson action
with equivariant momentum map 
$\check{L}\colon\check{P}\rightarrow G^{\ast}$
given by 
$$\check{L}(p,d)=J_{l}(d).\eqn\Lcheck$$
Furthermore, the following identities are valid: 

(1) $J_{r}\comp l_{k}=J_{r},\;\forall k\an G$;

(2) $J_{l}\comp r_{h}=J_{l},\;\forall h\an H$;

(3) $r_{h}\comp l_{k}=l_{k}\comp r_{h},\;\forall k\an G,h\an H$;

(4) $\check{l}_{k}\comp\check{\sigma}_{h}=\check{\sigma}_{h}\comp
\check{l}_{k},\;\forall k\an G,h\an H$.}

\undertext{\it Proof}. The fact that $\check{l}$ is a Poisson action
is evident. In order to prove that $\check{L}$ 
defined in $\Lcheck$ is a momentum map for $\check{l}$, it is
sufficient to apply Proposition {\lrpoissonmoment}(2)
choosing $\sigma_{1}$ as the trivial action of $G$ on $P$ and 
$\sigma_{2}=l$. In that case, the constant map $P\rightarrow 
G^{\ast}$ which to each point associates the identity of $G^{\ast}$, 
is an equivariant momentum map for $\sigma_{1}$.

Let now $k\an G$, $h\an H$ and $d=gu=u_{1}g_{1}\an D(G)$. Then:
$$\eqalign{(J_{r}\comp l_{k})(d)&=J_{r}(\lambda_{\rho_{g^{-1}}(u)}
(k)g\cdot u)\cr
\hfill&=s(i^{\ast}(u))\cr
\hfill&=J_{r}(d).\cr}$$
We check now relation (2):
$$\eqalign{(J_{l}\comp r_{h})(d)&=J_{l}(dh)\cr
\hfill&=J_{l}(u_{1}g_{1}h)\cr
\hfill&=u_{1}\cr
\hfill&=J_{l}(d),\cr}$$
because $J_{l}$ coincides with the projection 
$p_{2}^{t}\colon D(G)\rightarrow G^{\ast}$ defined by 
$p_{2}^{t}(u_{1}g_{1})=u_{1}$ (see {\jhlu}). We omit the proof of
(3) which is based on similar techniques. We finally check the 
validity of (4) making use of the commutativity between 
$r_{h}$ and $l_{k}$.\QED

\vskip0.3cm

We now observe that the momentum map $\check{J}$ given by
$\Jcheck$ is invariant under the action $\check{l}$: 
$(\check{J}\comp\check{l}_{k})(p,d)=J(p)J_{r}(l_{k}(d))=
\check{J}(p,d)$, $\forall k\an G,(p,d)\an\check{P}$, thanks to
relation (1) of Proposition {\invar}. Thus, we obtain an action
$\check{l}\colon G\times\check{J}^{-1}(e^{\ast})\rightarrow
\check{J}^{-1}(e^{\ast})$ which commutes with the 
action $\check{\sigma}$ of $H$ on the submanifold 
$\check{J}^{-1}(e^{\ast})$ (we recall that 
$\check{J}$ is equivariant, so we have an action
$\check{\sigma}\colon H\times\check{J}^{-1}(e^{\ast})\rightarrow
\check{J}^{-1}(e^{\ast})$). Consequently, we have a left action
$l_{ind}\colon G\times P_{ind}\rightarrow P_{ind}$ of $G$ on the
induced manifold $P_{ind}$. We will show that this action is Poisson.
To this end, it is more convenient to reformulate the Poisson 
property of an action in terms of differentiable functions. Thus,
using the Lie bracket on the 1-forms on a Poisson manifold 
$(P,\pi_{P})$ and the infinitesimal expression of the Poisson property 
of an action, one finds that the action
$\sigma\colon(G,\pi_{G})\times(P,\pi_{P})\rightarrow(P,\pi_{P})$
is Poisson if and only if 
$$\sigma(X)\{F,H\}=\{\sigma(X)F,H\}+\{F,\sigma(X)H\}
+(\sigma\wedget\sigma)\delta(X)(dF\otimes dH),\eqn\poissonvar$$
for each $X\an\frak g$, $F,H\an C^{\infty}(P)$,
where $\delta\colon\frak g
\rightarrow{\Lambda}^{2}\frak g$ is the 
linearization of $\pi_{G}$ at the identity of $G$. 
In our case, if $s\colon\check{J}^{-1}(e^{\ast})
\rightarrow P_{ind}$ is the projection, and 
$i_{e}\colon\check{J}^{-1}(e^{\ast})\hookrightarrow\check{P}$
the canonical inclusion, then, the following equation is valid
$$s^{\ast}\{F,H\}=i_{e}^{\ast}\{\tilde{F},\tilde{H}\},
\eqn\quopoisson$$
for each $F,H\an C^{\infty}(P_{ind})$, 
where $\tilde{F},\tilde{H}$
are arbitrary local extensions of 
$s^{\ast}F,s^{\ast}H$ respectively, such that 
$d\tilde{F},d\tilde{H}$ vanish on the 
subcharacteristic distribution 
${\cal C}\check{J}^{-1}(e^{\ast})$ (see {\mr}, {\va}).
Taking into account the fact that the infinitesimal generator
$l_{ind}(X)$ is obtained by projection of $\check{l}(X)$,
$\forall X\an\frak g$, one can write:
$$\eqalign{s^{\ast}(l_{ind}(X)\{F,H\})&=
\check{l}(X)(s^{\ast}\{F,H\})\cr
\hfill&=\check{l}(X)(i_{e}^{\ast}\{\tilde{F},\tilde{H}\})\cr
\hfill&=i_{e}^{\ast}(\{\check{l}(X)\tilde{F},\tilde{H}\}+
\{\tilde{F},\check{l}(X)\tilde{H}\}+(\check{l}
\wedget\check{l})\delta(X)
(d\tilde{F}\otimes d\tilde{H}))\cr
\hfill&=s^{\ast}(\{l_{ind}(X)F,H\}+\{F,l_{ind}(X)H\}\cr
\hfill&\kern11pt+(l_{ind}\wedget l_{ind})\delta(X)
(d\tilde{F}\otimes d\tilde{H})),\cr}$$
which confirms our assertion. Note that we used the 
fact that the function $\check{l}(X)\tilde{F}$ 
is an extension of $\check{l}(X)s^{\ast}F=s^{\ast}(l_{ind}(X)F)$
whose differential vanishes on vector fields taking their 
values in ${\cal C}\check{J}^{-1}
(e^{\ast})$: $d(\check{l}(X)\tilde{F})(\check{\sigma}(Y))=
\check{\sigma}(Y)\check{l}(X)\tilde{F}=\check{l}(X)\check{\sigma}(Y)
\tilde{F}=0$, thanks to the commutativity between 
the actions $\check{l}$ et $\check{\sigma}$ 
(Proposition {\invar}(4)) and to the fact that 
$d\tilde{F}$ vanishes on ${\cal C}\check{J}^{-1}(e^{\ast})$.

Consider now the momentum mapping  $\check{L}\colon
\check{P}\rightarrow G^{\ast}$ given by $\Lcheck$. By 
the invariance of the momentum mapping $J_{l}$ under the action
$r$, we easily find that $\check{L}$ is invariant under the action
$\check{\sigma}$. Thus, $\check{L}$ projects to a well-defined 
differentiable map $J_{ind}\colon P_{ind}\rightarrow G^{\ast}$.
Then, the defining equation $\check{l}(X)=\pi_{\check{P}}
^{\sharp}(\check{L}^{\ast}X^{l})$ of $\check{L}$, shows clearly
that we also have $l_{ind}(X)=\pi_{ind}^{\sharp}
(J_{ind}^{\ast}X^{l})$, $\forall X\an\frak g$, which means
that $J_{ind}$ is an equivariant momentum map for $l_{ind}$.
We have proved:

\math{Theorem.}{\proper}{\sl
Let $(G,\pi_{G})$ be a 
Poisson-Lie group, $(H,\pi_{H})$
a closed Poisson-Lie subgroup of $(G,\pi_{G})$ and $\sigma\colon
(H,\pi_{H})\times(P,\pi_{P})\rightarrow(P,\pi_{P})$ a proper left 
Poisson
action on the Poisson manifold $(P,\pi_{P})$, admiting the  
equivariant momentum map $J\colon P\rightarrow H^{\ast}$. If 
$(G,\pi_{G})$
is complete, connected and simply connected, then
there exists a Poisson manifold $(P_{ind},\pi_{ind})$, obtained
in $\Pind$ by reduction through the momentum mapping given by
$\Jcheck$, and a left Poisson action $l_{ind}\colon(G,\pi_{G})
\times(P_{ind},\pi_{ind})\rightarrow(P_{ind},\pi_{ind})$, induced
by the action $\lcheckaction$, admiting an equivariant 
momentum map given by $\Lcheck$. The manifold 
$(P_{ind},\pi_{ind})$ is called induced Poisson manifold.}

\vskip0.3cm

\goodbreak

{\bf Examples} 

{\bf (1) Poisson induction from a point.} We consider the case 
where the Poisson manifold $(P,\pi_{P})$ is a point with the zero 
Poisson structure: $P=\{\hbox{point}\}$ and the Poisson action of 
$(H,\pi_{H})$ is trivial with the momentum mapping $J\colon 
P\rightarrow H^{\ast}$ given by a fixed element $u_{0}\an H^{\ast}$: 
$J(p)=u_{0}$. The equivariance condition for such a momentum mapping 
is equivalent to the invariance of $u_{0}$ under the left dressing 
transformations of $H$ on $H^{\ast}$. Choosing now a Lie group morphism
$s^{\ast}\colon H^{\ast}\rightarrow G^{\ast}$
which commutes with left dressing 
transformations, we obtain a map 
$j=s^{\ast}\comp J\colon P\rightarrow G^{\ast}$ and let $j(p)=w_{0}$.
This defines, according to {\phd}, a diffeomorphism $I\colon
\check{J}^{-1}(e^{\ast})\rightarrow P\times G\times H^{0}$ given by 
$I(p,gu)=(p,guw_{0}^{-1})$, where $H^{\circ}\subset G^{\ast}$ is the fibre 
over the identity of the canonical projection $G^{\ast}\rightarrow 
H^{\ast}$. Under this identification, the induced Poisson manifold is 
diffeomorphic to the associated bundle $G\times_{\kern-3pt H\kern2pt}
H^{\circ}$, which 
carries a natural symplectic structure obtained either by Poisson 
reduction of the symplectic manifold $D(G)$ {\phd} or by the 
construction of the symplectic groupoid of the reduced Poisson space 
$G/H$ {\xu}. The Poisson induction procedure modifies this natural structure 
in the following manner. If $Q\colon D(G)\rightarrow D(G)$ is the 
diffeomorphism given by right multiplication with the element 
$w_{0}^{-1}$, then a direct calculation shows that 
$Q(\pi_{+}(d))=\pi_{+}(Q(d))+L_{d}\pi_{-}(w_{0}^{-1})$, which means 
that the Poisson induction from a point, leads to a modification of 
the canonical symplectic structure of the symplectic groupoid of $G/H$ 
(or, using the terminology of {\phd}, Poisson cotangent bundle of 
$G/H$) identified with the associated bundle 
$G\times_{\kern-3pt H\kern2pt}
H^{\circ}$. Clearly, the modification term vanishes when $u_{0}=e^{\ast}$, 
because $\pi_{-}$ is a Poisson-Lie structure on $D(G)$. 
This is the exact Poisson analog of the modified cotangent 
bundle of a homogeneous space $G/H$ {\det}.

\vskip0.5cm

{\bf (2) Poisson induced orbits.} We are placed now in the case where 
$P=H\cdot v$, the orbit of the element $v\an H^{\ast}$ under the right 
dressing transformations of $H$ on $H^{\ast}$. This action is 
Hamiltonian with momentum mapping given by the inclusion of $P$ in 
$H^{\ast}$. Let $w\an G^{\ast}$ be an element of the dual group of 
$G$ such that $i^{\ast}w=v$. We make the assumption that the fibre 
of $i^{\ast}\colon G^{\ast}\rightarrow H^{\ast}$ over $v$ is 
contained in the orbit of $w$ under the dressing action of  the 
subgroup $H$:
$$wH^{\circ}\subset H\cdot w.\eqn\orbitcondition$$
The constraint submanifold $\check{J}^{-1}(e^{\ast})$ consists in pairs
$(p,gu)\an P\times D(G)$ for which $p=i^{\ast}u$ and the action of 
$H$ on $\check{J}^{-1}(e^{\ast})$ is given by
$$\check{\sigma}_{h}(p,gu)=(\rho_{h^{-1}}(p),gh^{-1}\rho_{h^{-1}}(u)),$$
and therefore the equivalence class $[p,gu]$ must be written 
as $[p,gu]=\rho_{g^{-1}}(u)$. But if we write $p=\rho_{k^{-1}}(v)$, 
$k\an H$,
then $u=\rho_{k^{-1}}(wu^{\circ})$, $u^{\circ}\an H^{\circ}$, because 
$H^{\circ}$ is invariant under the dressing action of $H$. Taking 
into account the condition $\orbitcondition$, we find
$$P_{ind}=G\cdot w,$$
that is the orbit of $w\an G^{\ast}$ under the right dressing 
transformations of $G$, is obtained by Poisson induction on the orbit 
of $v=i^{\ast}w\an H^{\ast}$ under the dressing transformations of 
$H$. 

{\bf Remark.} {\sl The previous example is a Poisson generalization of 
one of the main results 
of {\rm\bag} concerning the geometry of the coadjoint orbits of a 
semi-direct product. Indeed, it is shown in {\rm\bag} that each 
coadjoint orbit of a semi-direct product can be obtained by 
symplectic induction on a coadjoint orbit of a conveniently choosen 
subgroup. The symplectic construction, concerning semi-direct 
products, 
is obtained from the Poisson one discussed here, 
if one takes (in the notation of {\rm\bag})
$G=K\times_{\kern-3pt 
\rho\kern2pt}V$ and $H=K_{p}\times_{\kern-3pt 
\rho\kern2pt}V$, both with the zero Poisson structure. Let us note 
that for a semi-direct product, the condition $\orbitcondition$ is 
satisfied for all the elements $w$.}

\vfill\eject
\vskip1cm
   \ifreferenceopen \Closeout\referencewrite \referenceopenfalse \fi
   \line{\bf\hskip0pt\hfil References\hfil}\vskip\headskip
   \vskip0.3cm
   \input referenc.txa

\end